\documentclass[10pt]{amsart}

\newtheorem{thm}{Theorem}[section]
\newtheorem{prop}{Proposition}[section]

\newtheorem{obs}{Remark}[section]
\newtheorem{defin}{Definition}[section]

\usepackage{amssymb}
\usepackage{amsmath}
\usepackage{amsfonts}
\usepackage{graphicx}

\numberwithin{equation}{section}

\usepackage[T1]{fontenc}
\usepackage{palatino}

\begin{document}
\title[geometric measure and 3D NSE]
{A geometric measure-type regularity criterion for solutions to
the 3D Navier-Stokes equations}
\author{Z. Gruji\'c}
\address{Department of Mathematics\\
University of Virginia\\ Charlottesville, VA 22904}
\date{\today}
\begin{abstract}
A local anisotropic geometric measure-type condition on the super-level sets of
solutions to the 3D NSE preventing the formation of a finite-time
singularity is presented; essentially, \emph{local
one-dimensional sparseness} of the regions of intense fluid activity
in a very weak sense.

\end{abstract}

\maketitle

\newpage

\section{Prologue}

\noindent Rigorous study of geometric depletion of the
nonlinearity in the 3D Navier-Stokes equations (3D NSE) was
initiated by Constantin in \cite{Co94}; the approach was based on
the singular integral representation formula for the stretching
factor in the evolution of the vorticity magnitude featuring a
geometric kernel depleted by coherence of the vorticity direction.
This representation was subsequently utilized by Constantin and
Fefferman in \cite{CoFe93} to show that as long as the vorticity
direction is Lipschitz-coherent, no finite-time blow up can occur
and later by Beirao da Veiga and Berselli in \cite{daVeigaBe02},
where the Lipschitz-coherence regularity condition was replaced with
$\frac{1}{2}$-H\"older.

\medskip

Spatiotemporal localization of the $\frac{1}{2}$-H\"older-coherence
regularity
criterion was performed in \cite{GrZh06, Gr09}, and also -- utilizing
a different
localization procedure -- by Chae, Kang and Lee in \cite{ChKaLe07}.

\medskip

The $\frac{1}{2}$-H\"older-coherence condition is super-critical with respect
to the natural scaling of the 3D NSE; a family of scaling-invariant, critical,
hybrid
geometric-analytic local regularity criteria -- including a scaling-invariant
improvement of the $\frac{1}{2}$-H\"older-coherence condition -- was presented
in \cite{GrGu10-1}.

\medskip

In the realm of the mathematical theory of turbulence, the
$\frac{1}{2}$-H\"older-coherence
condition was recently (\cite{DaGr11-1}) paired with the condition
on a modified Kraichnan scale
to obtain
a first rigorous evidence of existence of (anisotropic) enstrophy
cascade in 3D viscous incompressible flows.

\medskip

A different approach to discovering geometric scenarios ruling out
formation of singularities in the 3D NSE was introduced in
\cite{Gr01}. The main idea was to utilize the local-in-time spatial
analyticity properties of solutions in $L^p$ (\cite{GrKu98}) via the
plurisubharmonic measure maximum principle -- a generalization to
$\mathbb{C}^n$ (cf. \cite{Sad81}) of the classical harmonic measure
majorization principle in the complex plane (the log-convexity of
the modulus of an analytic function; see, e.g., \cite{Nev70}).

\medskip

The regularity criterion derived in \cite{Gr01} is a condition on
the regions of intense fluid activity near a possible blow up time
requiring local existence of a sparse coordinate \emph{projection}
on the scale comparable to the uniform radius of spatial
analyticity. The estimate on the plurisubharmonic measure was
performed within the framework of product-type domains -- hence the
requirement on a coordinate projection. This could be somewhat
relaxed, but not substantially due to the rigidity of the
$\mathbb{C}^n$ structure. Also, once the computation of the
plurisubharmonic measure was reduced to the computation of the
individual (coordinate) harmonic measures, the estimate on the
harmonic measure was carried out with respect to an infinite strip,
giving the argument a nonlocal character.

\medskip

In the present work, we completely bypass the rigidity of the
$\mathbb{C}^n$ structure, resulting in a much weaker \emph{local
geometric measure-type condition}. Utilizing translational and
rotational invariance of the 3D NSE, as well as some basic geometric
properties of the harmonic measure, the argument is ultimately
reduced to the problem of estimating the harmonic measure of an
arbitrary closed subset of $[-1,1]$ computed at $0$ with respect to
the unit disk. This is a generalization of the classical Beurling's
problem \cite{Beu33, Nev70} proposed by Segawa in \cite{Seg88}; a
symmetric version was solved by Essen and Haliste in \cite{EssHa89},
and the general case relatively recently by Solynin in \cite{Sol99}
via a general symmetrization argument. The utility of Solynin's
result here is that it allows us to formulate the condition in view
in terms of the ratio of one-dimensional Lebesgue measures, rather
than in terms of a specific `rearrangement' of the one-dimensional
trace of the region of intense fluid activity. The engine
behind the proof is the interplay between the diffusion in the model
-- quantified by the local-in-time (sharp) analytic smoothing in
$L^\infty$ -- and the geometric properties of the harmonic measure.

\medskip

A precise statement of our regularity criterion will be given in the
main text; at this point, we convey the essence of the result.
Denoting the region of intense fluid activity at time $s$ by
$\Omega_s(M)$ and the lower bound on the uniform radius of spatial
analyticity by $\rho(s)$ (this take places near a potential singular
time), the condition in view is simply a stipulation that for a
given point $x_0$, there exists a radius $r=r(x_0)$, $0 < r \le
\rho(s)$ and a unit vector $d=d(x_0)$, such that

\begin{equation}\label{yay}
  \frac{|\Omega_s(M) \cap (x_0-rd, x_0+rd)|}{{2r}} \le
  \delta
\end{equation}
for some $\delta$ in $(0,1)$. (It will transpire that it is enough
to require (\ref{yay}) on a suitably chosen finite sequence of times
.) There are two versions of the result, one for the velocity and
one for the vorticity, based on the spatial analyticity estimates on
solutions to the velocity and vorticity formulations of the 3D NSE,
respectively.

\medskip

It is plain that (\ref{yay}) is a much weaker condition than the one
in \cite{Gr01}; all that is needed here is local sparseness of a
one-dimensional \emph{trace} of the region of intense fluid activity
in a very weak sense. On the other hand, it is of a different nature
and hence not directly comparable to the coherence of the vorticity
direction-type regularity criteria.

\medskip

The main efficacy of the aforementioned regularity criterion is in
ruling out various sparse geometric scenarios for a finite-time blow
up, both in the velocity and the vorticity formulations. As
succinctly put by P. Constantin, ``intermittency implies
regularity'' \cite{PC11}. An example of interest that can be ruled
out is a blow up scenario in which the region of intense vorticity
(defined as a region in which the vorticity magnitude exceeds a
suitable fraction of the $L^\infty$-norm) is -- at suitable
near-blow up times -- comprised of vortex filaments with diameters
of the cross-sections bounded above by the uniform radius of spatial
analyticity.

\medskip

The paper is organized as follows. In Section 2, we collect relevant
properties of the harmonic measure in the plane, and in Section 3 we
recall the local-in-time spatial analyticity of solutions in
$L^\infty$. Section 4 contains the main result. The last section
indicates a scenario which -- in a statistically significant sense
-- leads to closing the scaling gap in the regularity problem.

\medskip

\section{Harmonic measure}

Basic properties of the harmonic measure in the complex plane can be
found, e.g., in \cite{Nev70, Ahl10}. First, we briefly recall a few
relevant facts following \cite{Ahl10}.

\medskip

Let $\Omega$ and $K$ be an open and a closed set in the complex
plane, respectively. When the geometry of $\Omega \setminus K$ is
not too convoluted, there exists a unique bounded harmonic function
on $\Omega \setminus K$, denoted by $\omega=\omega(\cdot,\Omega,K)$,
such that -- in the sense of a well-defined limit as a point
approaches the boundary -- $\omega$ is equal to 1 on $K$ and 0 on
the rest of the boundary; $\omega(z,\Omega,K)$ is the harmonic
measure of $K$ with respect to $\Omega$ computed at $z$.

\medskip

Two straightforward consequences of the general harmonic measure
majorization principle, c.f. Theorem 3.1 \cite{Ahl10}, are the
following (\cite{Ahl10}, p. 39).

\begin{prop}\label{mono}
The harmonic measure is increasing (as a measure) with respect to
both $K$ and $\Omega$.
\end{prop}

\begin{prop}\label{2_c}
Let $f$ be analytic in $\Omega \setminus K$, $|f| \le M$, and $|f|
\le m$ on $K$ (in the sense of \, $\limsup$ as a point
approaches the boundary). Then
\[
 |f(z)| \le m^\theta \, M^{1-\theta}
\]
for any $z$ in $\Omega \setminus K$, where
$\theta=\omega(z,\Omega,K)$.
\end{prop}

This a refined form of the maximum modulus principle for analytic
functions in $\Omega \setminus K$ (the log-convexity of the modulus
of $f$ -- sometimes referred to as ``two-constants theorem'').

\medskip

Another useful property of the harmonic measure is the following
(see., eg., \cite{Nev70}).

\begin{prop}\label{conf}
The harmonic measure is invariant with respect to conformal
mappings.
\end{prop}

Finally, we recall a result on extremal properties of the harmonic
measure in the unit disk $\mathbb{D}$ obtained by Solynin in
\cite{Sol99}.

\begin{thm}\label{K}
Let $K$ be a closed subset of $[-1,1]$ such that $|K| = 2\lambda$
for some $\lambda$, $0<\lambda<1$,
and suppose that $0$ is in $\mathbb{D} \setminus K$. Then
\[
 \omega(0,\mathbb{D},K) \ge \omega(0,\mathbb{D}, K_\lambda) =
 \frac{2}{\pi} \arcsin \frac{1-(1-\lambda)^2}{1+(1-\lambda)^2}
\]
where $K_\lambda = [-1, -1+\lambda] \cup [1-\lambda, 1]$.
\end{thm}

The above theorem provides a generalization of the classical
Beurling's result \cite{Beu33} in which $K$ is a finite union of
intervals lying on one side of the origin. This was conjectured by
Segawa in \cite{Seg88}, and the symmetric version was previously
resolved in \cite{EssHa89}.

\medskip

\section{Spatial analyticity in $L^\infty$}

The 3D NSE equations read
\begin{equation}\label{nse_u}
u_t + (u \cdot \nabla)u = -\nabla p + \triangle u
\end{equation}
supplemented with the incompressibility condition $\nabla \cdot u = 0$,
where $u$ is the velocity of the fluid and $p$ the pressure (the
viscosity is set to 1).

\medskip

A method for deriving explicit local-in-time lower bounds on the
uniform radius of spatial analyticity of solutions to the NSE in
$L^p$ was introduced in \cite{GrKu98}; see also \cite{Ku99} for
analogous results in the vorticity formulation. We will make use of
the following sharp analyticity estimate in $L^\infty$ (cf.
\cite{Gu10}; \cite{Ku03} for the corresponding real result).

\begin{thm}\label{an_u}
Let $u_0$ be in $L^\infty(\mathbb{R}^3)$. Then, there exists an
absolute constant $c_0 > 1$ such that setting
$\displaystyle{T=\frac{1}{c_0^2 \|u_0\|_\infty^2}}$, a unique mild
solution $u=u(t)$ on $[0,T]$ has the analytic extension $U=U(t)$ to
the region
\[
 \mathcal{R}_t = \{x+iy \in \mathbb{C}^3 : \, |y| \le \frac{1}{c_0}
 \sqrt{t}\}
\]
for any $t$ in $(0,T]$. In addition,
\[
 \|U(t)\|_{L^\infty(\mathcal{R}_t)} \le c_0 \|u_0\|_\infty
\]
for all $t$ in $[0,T]$.
\end{thm}

\begin{obs}
\emph{Recall that -- given a divergence-free initial datum $u_0$ in
$L^\infty$ -- a $C\bigl((0,T); L^\infty\bigr)$ function $u$ is a
mild solution to the 3D NSE corresponding to $u_0$ provided it
solves
\begin{align*}
 u^k(x,t) = & - \iint \partial_j K (x-y,t-s) \, u^j(y,s) u^k(y,s) \, dy \, ds\\
                & -\iint \partial_k K(x-y,t-s) \, p(y,s) \, dy \, ds\\
                & + \int K(x-y,t) \, u_0^k(y) \, dy,
\end{align*}
where $K$ is the heat kernel and $p=-R_iR_j u^i u^j$ ($R_k$ being
the $k$-th Riesz Transform).}
\end{obs}

The vorticity formulation of the 3D NSE reads
\begin{equation}\label{nse_omega}
\omega_t + (u \cdot \nabla)\omega = (\omega \cdot \nabla)u +
\triangle \omega
\end{equation}
where $\omega = \, \mbox{curl} \, u$ is the vorticity; the vorticity
version of the above theorem is as follows (the proof
is analogous; utilizing the Biot-Savart law to close each
iteration).

\medskip

\begin{thm}\label{an_omega}
Let $\omega_0$ be in $L^\infty(\mathbb{R}^3)$. Then, there exists an
absolute constant $d_0 > 1$ such that setting
$\displaystyle{T=\frac{1}{d_0^2 \|\omega_0\|_\infty}}$, a unique
mild solution $\omega=\omega(t)$ on $[0,T]$ has the analytic
extension $\Omega=\Omega(t)$ to the region
\[
 \mathcal{R}_t = \{x+iy \in \mathbb{C}^3 : \, |y| \le \frac{1}{d_0}
 \sqrt{t}\}
\]
for any $t$ in $(0,T]$. In addition,
\[
 \|\Omega(t)\|_{L^\infty(\mathcal{R}_t)} \le d_0 \|\omega_0\|_\infty
\]
for all $t$ in $[0,T]$.
\end{thm}

\begin{obs}\label{abs}
\emph{ An inspection of the proofs in \cite{Gu10, Ku03} reveals that
$c_0$ (and similarly $d_0$) is an \emph{absolute} constant depending
only on the $BMO$ bound on the Riesz Transforms in $\mathbb{R}^3$,
and the bound appearing in the result on the non-homogeneous heat
equation establishing $L^\infty$-regularity in the case the
non-homogeneity is a divergence of a $BMO$ function (see, e.g.,
Lemma 3.1 in \cite{Ku03}).}
\end{obs}

\section{The main result}

We start with introducing a geometric measure-theoretic concept of
\emph{weak local linear sparseness of a set around a point, at a
given scale}, suitable for our purposes.

\medskip

\begin{defin}
Let $x_0$ be a point in $\mathbb{R}^3$, $r>0$, $S$ an open subset of
$\mathbb{R}^3$ and $\delta$ in $(0,1)$.

\medskip

The set $S$ is \emph{linearly $\delta$-sparse around $x_0$ at scale
$r$ in weak sense} if there exists a unit vector $d$ in $S^2$ such
that
\[
 \frac{|S \cap (x_0-rd, x_0+rd)|}{2r} \le \delta.
\]
\end{defin}

\medskip

In what follows, we derive the main result for the velocity
formulation and simply state the analogous result for the vorticity
formulation; modifying the proof in the second case is essentially
relabeling.

\medskip

For $M>0$, denote by $\Omega_t(M)$ the super-level set at time $t$;
more precisely,
\[
 \Omega_t(M) = \{x \in \mathbb{R}^3: |u(x,t)| > M\}.
\]
Then, our main result reads as follows.

\medskip

\begin{thm}\label{sparse_u}
Suppose that a solution $u$ is regular on an interval $(0,T^*)$.
(Recall that $u$ is then necessarily in $C\Bigl((0,T^*);
L^\infty\Bigr)$.)

\medskip

Assume that either

\medskip

(i) \ there exists $t$ in $(0,T^*)$ such that
$\displaystyle{t+\frac{1}{c_0^2 \|u(t)\|_\infty^2} \ge T^*}$
($c_0$ is the constant
featured in Theorem \ref{an_u}), or

\medskip

(ii) \ $\displaystyle{t+\frac{1}{c_0^2 \|u(t)\|_\infty^2} < T^*}$
for all $t$ in $(0,T^*)$,
and there exists $\epsilon$ in $(0,T^*)$ such that for any
$t$ in $(T^*-\epsilon, T^*)$, there
exists $s=s(t)$ in $\Bigl[t+\frac{1}{4c_0^2 \|u(t)\|_\infty^2},
t+\frac{1}{c_0^2 \|u(t)\|_\infty^2}\Bigr]$
with the property that for any spatial point
$x_0$, there exists a scale $r=r(x_0)$, $0<r\le \frac{1}{2c_0^2
\|u(t)\|_\infty}$, such that the super-level set
$\Omega_s(M)$ is linearly $\delta$-sparse around $x_0$ at scale $r$
in weak sense; here, $\delta=\delta(x_0)$ is an arbitrary value
in $(0,1)$,
$h=h(\delta)=\frac{2}{\pi}\arcsin\frac{1-\delta^2}{1+\delta^2}$,
$\alpha=\alpha(\delta)\ge\frac{1-h}{h}$, and
$M=M(\delta)=\frac{1}{c_0^\alpha}
\|u(t)\|_\infty$.

\medskip

Then, there exists $\gamma >0$ such that $u$ is in
$L^\infty\Bigl((T^*-\epsilon, T^*+\gamma); L^\infty\Bigr)$, i.e.,
$T^*$ is not a singular time.
\end{thm}

\begin{proof} There are two cases to consider.

\medskip

\noindent Case (i)

\medskip

In this case, the statement of the theorem follows from
Theorem \ref{an_u}, setting the initial time to $t$.

\medskip

\noindent Case (ii)

\medskip

Pick a time $t_0$ in $(T^*-\epsilon, T^*)$, and let $s_0=s(t_0)$ be as in
the statement of the theorem. Then, for
any $x_0$ in $\mathbb{R}^3$, there exists $r=r(x_0)$ , $0 < r \le
\frac{1}{2c_0^2 \|u(t_0)\|_\infty}$, a direction vector
$d=d(x_0)$, and $\delta = \delta(x_0)$ in $(0,1)$ such that

\medskip

\[
 \frac{|\Omega_{s_0}(M) \cap (x_0-rd, x_0+rd)|}{2r} \le
 \delta
\]
($M=M(\delta)$ as in the statement).

\medskip

The intent is to show
\begin{equation}\label{goal}
 \|u(s_0)\|_\infty \le \|u(t_0)\|_\infty.
\end{equation}

\medskip

Fix $x_0$. Recall that the NSE exhibit translational and rotational
invariance in the spatial variable. Translate for $-x_0$, rotate by
the matrix $Q$ transforming the unit direction $d$ to the coordinate
vector $e_1$ and denote the transformed solution by $u_{x_0,Q}$.
Then, $u_{x_0,Q}(x,t)=Q u\bigl(Q^{-1}(x+x_0),t\bigr)$.

\medskip

Solve the NSE locally-in-time starting at $t_0$; the spatial
analyticity properties of the solution at time $s_0$ are given by
simply translating in time the statement of Theorem \ref{an_u}.

\medskip

Moreover, since the rotation $Q$ has no effect on computing the
norms, the transformed solution $u_{x_0,Q}$ at time $s_0$ enjoys
exactly the same analyticity features.

\medskip

In particular -- focusing on the first coordinate --
$u_{x_0,Q}(s_0)$ is spatially analytic on a strip, symmetric around
the real axis, with the width equal to (at least)
\[
 \rho(s_0)= \frac{1}{2 c_0^2 \|u(t_0)\|_\infty}.
\]
The region of interest is the disk around the origin with the radius
$r$, $D_r$. Note that $D_r$ is contained in the domain of
analyticity of $u_{x_0,Q}(s_0)$.

\medskip

Our goal is to obtain an improved estimate on $u_{x_0,Q}(0,s_0)$.
Denote by $K$ the complement of the image of the set
$\Omega_{s_0}(M) \cap (x_0-rd, x_0+rd)$, under the change of
coordinates, in $[-r,r]$. Then, $K$ is closed, and the sparseness
assumption implies $|K| \ge 2r(1-\delta)$. If $0$ is in $K$,
$|u_{x_0,Q}(0,s_0)| < \|u(t_0)\|_\infty$, and we are done (with this
$x_0$). If not, the harmonic measure maximum principle -- Proposition
\ref{2_c} -- together with the $L^\infty$-bound on the complexified
solution stated in Theorem \ref{an_u}, implies

\medskip

\begin{equation}\label{b1}
 |u_{x_0,Q}(0,s_0)| \le \Bigl(\frac{1}{c_0^\alpha}
 \|u(t_0)\|_\infty\Bigr)^{\omega(0,D_r,K)} \Bigl(c_0
 \|u(t_0)\|_\infty\Bigr)^{1-\omega(0,D_r,K)}.
\end{equation}

\medskip

Recall that the harmonic measure is invariant under conformal
mappings (Proposition \ref{conf}). In particular, it is invariant
under the scaling map $z \mapsto \frac{1}{r} z$. This paired with
the monotonicity of the harmonic measure with respect to $K$
(Proposition \ref{mono}) and Theorem \ref{K} yields

\medskip

\begin{equation}\label{stuff}
 \omega(0,D_r,K) \ge
 \frac{2}{\pi}\arcsin\frac{1-\delta^2}{1+\delta^2} = h.
\end{equation}

\medskip

Combining the estimates (\ref{b1}) and (\ref{stuff}) leads to

\medskip

\begin{equation}\label{b2}
|u_{x_0,Q}(0,s_0)| \le \Bigl(\frac{1}{c_0^\alpha}
 \|u(t_0)\|_\infty\Bigr)^h \Bigl(c_0
 \|u(t_0)\|_\infty\Bigr)^{1-h} \le \|u(t_0)\|_\infty.
\end{equation}

\medskip

This, in turn, implies $|u(x_0,s_0)| \le \|u(t_0)\|_\infty$, and
since $x_0$ was an arbitrary spatial point in $\mathbb{R}^3$,
$\|u(s_0)\|_\infty \le \|u(t_0)\|_\infty$.

\medskip

Let
\[
 M_0=\|u(t_0)\|_\infty.
\]

Setting $t_1=s_0$ and repeating the argument yields $\|u(s_1)\|_\infty
\le \|u(t_1)\|_\infty \le M_0$, where $s_1=s(t_1)$. After finitely many steps,
we reach the time $s_n, s_n < T^*$ such that
$\displaystyle{\|u(s_n)\|_\infty \le M_0 \ \mbox{and} \ s_n +
\frac{1}{c_0^2 M_0^2} > T^*}$.
The statement of the theorem now follows from
Theorem \ref{an_u}, setting the initial time to $s_n$.
\end{proof}

\medskip

\begin{obs}
\emph{It is plain from the proof that it is enough to assume the
condition on a \emph{finitely many} suitably chosen times.}
\end{obs}

\medskip

The vorticity version and the proof are completely analogous --
utilizing Theorem \ref{an_omega} in place of Theorem \ref{an_u}.

\medskip

For $M>0$, denote by $\Omega^\omega_t(M)$ the vorticity super-level
set at time $t$; more precisely,
\[
 \Omega^\omega_t(M) = \{x \in \mathbb{R}^3: |\omega(x,t)| > M\}.
\]

\medskip

\begin{thm}\label{sparse_omega}
Suppose that a solution $u$ is regular on an interval $(0,T^*)$.

\medskip

Assume that either

\medskip

(i) \ there exists $t$ in $(0,T^*)$ such that
$\displaystyle{t+\frac{1}{d_0^2 \|\omega(t)\|} \ge T^*}$
($d_0$ is the constant
featured in Theorem \ref{an_omega}), or

\medskip

(ii) \ $\displaystyle{t+\frac{1}{d_0^2 \|\omega(t)\|} < T^*}$
for all $t$ in $(0,T^*)$,
and there exists $\epsilon$ in $(0,T^*)$ such that for any
$t$ in $(T^*-\epsilon, T^*)$, there
exists $s=s(t)$ in $\Bigl[t+\frac{1}{4d_0^2 \|\omega(t)\|},
t+\frac{1}{d_0^2 \|\omega(t)\|}\Bigr]$
with the property that for any spatial point
$x_0$, there exists a scale $r=r(x_0)$, $0<r\le \frac{1}{2d_0^2
\|\omega(t)\|_\infty^\frac{1}{2}}$, such that the super-level set
$\Omega^\omega_s(M)$ is linearly $\delta$-sparse around $x_0$ at scale $r$
in weak sense; here, $\delta=\delta(x_0)$ is an arbitrary value
in $(0,1)$,
$h=h(\delta)=\frac{2}{\pi}\arcsin\frac{1-\delta^2}{1+\delta^2}$,
$\alpha=\alpha(\delta)\ge\frac{1-h}{h}$, and
$M=M(\delta)=\frac{1}{d_0^\alpha}
\|\omega(t)\|_\infty$.

\medskip

Then, there exists $\gamma >0$ such that $\omega$ is in
$L^\infty\Bigl((T^*-\epsilon, T^*+\gamma); L^\infty\Bigr)$, i.e.,
$T^*$ is not a singular time.
\end{thm}

\medskip

\section{Epilogue}

Direct numerical simulations of turbulent flows reveal (see, e.g.,
\cite{SJO91}) that the preferred geometric signature of the regions
of intense vorticity is the one of vortex filaments. The general
agreement seems to be that the length of a filament is -- in a
statistically significant sense -- comparable with the macro scale.
For rigorous mathematical results
concerning creation and dynamics of vortex tubes in turbulent flows,
the reader is referred to \cite{CPS95}.

\medskip

Let us for a moment adopt the aforementioned geometry as a blow up
scenario. In order to make the reasoning more transparent, and avoid
stepping on the quicksand, set $\delta$ to be $\frac{1}{\sqrt{3}}$;
the corresponding values for $h$ and $\alpha$ are then $\frac{1}{3}$
and $2$, respectively. Recalling that $d_0$ is also an absolute
constant (cf. Remark \ref{abs}), all the constants appearing in
Theorem \ref{sparse_omega} are now absolute constants to be denoted
by either $C_i$ or $\frac{1}{C_i}$ for some $C_i > 1$. The
`constants' depending only on the initial data will be denoted by
$C_i^0$.

\medskip

Consider \emph{the region of intense vorticity} at a near-blow up
time $s(t)$ to be the region in which the vorticity magnitude
exceeds $\frac{1}{C_1} \|\omega(t)\|_\infty$. Then, Theorem
\ref{sparse_omega} implies that as long as the diameters of the
filaments' cross-sections are dominated by
$\displaystyle{\frac{1}{C_2}
\frac{1}{\|\omega(t)\|_\infty^{\frac{1}{2}}}}$, no blow up can
occur. At this point, recall that starting with the initial
vorticity a finite Radon measure, the $L^1$-norm of the vorticity is
bounded -- uniformly in time -- over any interval $(0,T)$)
\cite{Co90}. Tchebyshev inequality then implies the decrease of the
distribution function of the vorticity of at least
$\frac{C_3^0}{\lambda}$; consequently, the volume of the region of
intense vorticity decreases at least as $\displaystyle{C_4^0
\frac{1}{\|\omega(s(t))\|_\infty}}$. Assuming that the length of a
filament is comparable with the macro scale, this implies the
decrease of the diameters of the filaments' cross-sections of at
least $\displaystyle{C_5^0
\frac{1}{\|\omega(s(t))\|_\infty^{\frac{1}{2}}} \le C_5^0
\frac{1}{\|\omega(t)\|_\infty^{\frac{1}{2}}}}$, which is exactly
\emph{the scale} needed for the application of Theorem
\ref{sparse_omega} (without the loss of generality,
$\displaystyle{\|\omega(s(t))\|_\infty \ge \|\omega(t)\|_\infty}$;
if not, in the proof, one can simply take $s(t)$ as the new $t$).

\medskip

The above ruminations offer a geometric scenario leading to closing
the scaling gap in the regularity problem, i.e., to a manifestation
of \emph{criticality} for \emph{large data}. Assuming that the
`shape', i.e., the general geometry is correct, the weakest link is
the assumption that the length of a filament be comparable to the
macro scale; this was simply borrowed from the picture painted by
the numerical simulations. However, in a recent work
\cite{DaGr11-2}, the authors utilized a \emph{multiscale ensemble
averaging process} introduced in their study of turbulent cascades
in physical scales of 3D incompressible flows (\cite{DaGr10}) to
show that the averaged vortex-stretching term is -- near a possible
blow up time $T^*$ -- positive across a range of scales extending
from a power of a modified Kraichnan scale to the macro scale. This
provides a mathematical evidence of creation and persistence of the
macro scale-long vortex filaments (in a statistically significant
sense), and the pertaining research will be pursued in the future.

\medskip

\vspace{.2in}

\noindent ACKNOWLEDGMENTS The author thanks Professor Peter
Constantin for inspiring discussions and Department of Mathematics
at the University of Chicago for hospitality while being a Long Term
Visitor in Fall 2011; the support of the Research Council of Norway
via the grant number 213474/F20 and the National Science Foundation
via the grant number DMS-1212023 are gratefully acknowledged. In
addition, the author thanks the referees for their constructive
criticism.

\end{document}